\documentclass{amsart}

\title[Automorphism groups of hypergraphs]{Subgroups of simple primitive permutation groups defined by unordered relations }
\author{Mariusz Grech, Andrzej Kisielewicz}
\address{Institute of Mathematics, University of Wroclaw \\ 
pl.Grunwaldzki 2, 50-384 Wroclaw, Poland}
\email{Mariusz.Grech@math.uni.wroc.pl}
\thanks{Supported in part by Polish NCN grant 2016/21/B/ST1/03079; this is the full version of the extended abstract \emph{Automorphism groups of hypergraphs} presented in ICNAAM'19 \cite{GK2}}


\begin{document}


\newtheorem{Theorem}{Theorem}
\newtheorem{Lemma}{Lemma}
\newcommand{\G}{\mathcal G}


\maketitle


\begin{abstract}
The problem of describing the invariance groups of unordered relations, called briefly \emph{relation groups}, goes back to classical work by H. Wielandt. In general, the problem turned out to be hard, and so far it has been settled only for a few special classes of permutation groups. The problem has been solved, in particular, for the class of primitive permutation groups,  using the classification of finite simple groups and other deep results of permutation group theory. In this paper we show that, if $G$ is a finite simple primitive permutation group other than the alternating group $A_n$, then each subgroup of $G$, with four exceptions, is a relation group.
\end{abstract}

\section{INTRODUCTION}

A \emph{hypergraph} $\mathcal H$ is a family of subsets of a given set $\Omega$. The automorphism group of $\mathcal H$, denoted $Aut({\mathcal H})$, consists of permutations $g$ of $\Omega$ such that for every set $x\in\mathcal H$, $x^g\in \mathcal H$. The hypergraph $\mathcal H$ may be viewed as unordered relation on $\Omega$ and $Aut({\mathcal H})$ as the invariance group of this relation. That is why, permutation groups represented in this way are called just \emph{relation groups} (see \cite{SV}). If $\Omega$ is finite and we treat it as the set of indices $\Omega =\{1,2,\ldots,n\}$, then $\mathcal H$ defines a Boolean function  $f(x_1,x_2,\ldots,x_n)=1$ if and only if $\{i : x_i=1\} \in {\mathcal H}$. In this way, $Aut({\mathcal H})$ may be treated as the symmetry group of a Boolean function (see \cite{kis}). So there are three different ways to view relation groups, and we make use of various results obtained with different terminology. If $\mathcal H$ consists of two-element sets, then $\mathcal H$ is a graph and $Aut({\mathcal H})$ is the automorphism group of a graph.

A class of permutation groups for which it has been described exactly which of them can be represented as the automorphism groups of graphs or of other graphical structures, is the class of \emph{regular permutation groups}.  This is due to the one-to-one correspondence with the abstract groups and close connection with Cayley graphs. In fact, the problem in question has been solved for this class as a result of a large project under the name \emph{Graphical Regular Representation} of abstract groups. The final result, describing regular permutation groups that are the automorphism groups of graphs, has been published in Godsil \cite{god}. Based on this, 
similar results for directed graphs, tournaments, and oriented graphs have been obtained (see e.g., \cite{MS}). 
Partial results have been obtained also for the classes of cyclic permutation groups, abelian permutation groups, and Frobenius groups (see \cite{gre,GK1,spi}).

The question of representing permutation groups as automorphism groups of hypergraphs, that is, characterizing the relation groups, seems a little bit easier. It has been answered completely for regular groups and abelian groups \cite{GK}, and on the other extreme for primitive permutation groups in \cite{BC,SV}. The latter is based on advanced results of permutation group theory on primitive groups. 
The paper by Siemons and Volta contains also an effective tool \cite[Lemma~3.1]{SV} allowing to prove that every subgroup of a group $G$ in question is a relation group (providing $G$ satisfies some mild conditions). This gives hope for considerable progress towards the full characterization of the relation groups. 
In this paper we prove that, with four exceptions, each subgroup of a finite simple primitive permutation group other than $A_n$ is the automorphism group of a hypergraph.

\section{RESULT}
All permutation groups are finite and considered up to \emph{permutation isomorphism}. 
The symmetric and alternating groups on $n$ elements are denoted $S_n$ and $A_n$, respectively. The cyclic group generated by an $n$-element cycle is denoted $C_n$. We use standard notation, like $A(d,q)$ and $PSL(d,q)$, for affine and projective groups. Given a permutation group $G$, the set it acts on is denoted traditionally by $\Omega$. We use the notation $G^+$ for one point extension  obtained by adding an extra fix point to the group $G$.



Since we use tools worked out in \cite{SV}, generally, we follow the notation and terminology of \cite{SV}, and we refer the reader to it for a more detailed introduction to the topic. In particular, rather than about the automorphism groups $Aut({\mathcal H})$ of hypergraphs we speak of relation groups ${\G}(R)$ generated by a family $R$ of subsets of $\Omega$, called unordered relation (essentially, $R=\mathcal H$).

We should note that the groups defined by (ordered) relations have been introduced by Wielandt in \cite{wie}, and it must be emphasized that we are concerned  with \emph{unordered} relations and actions of permutation groups on (unordered) set. In this connection, recall that a permutation group $G < Sym(\Omega)$ is called \emph{set-transitive} if for every two subsets $x,y \subseteq \Omega$ of the same cardinality there exists a permutation $g\in G$ such that $xg = y$. This means that $G$ has exactly $n=|\Omega|$ orbits on the nonempty subsets of $\Omega$.  By the result of \cite{BP}, if $n=|\Omega|$ is finite, then apart from the symmetric and alternating groups, there are only 4 permutation groups (of degree $n=5,6$ or $9$) that are set-transitive.

Our result below is based on extensive knowledge we have on primitive groups. In particular, a crucial role is played by 
regular sets. Given a permutation group $G \leq Sym(\Omega)$, a subset $x\subseteq \Omega$ is called \emph{regular}, if the setwise stabilizer $G_x$ of the set $x$ in $G$ is trivial. In such a case $G$ acts regularly on the orbit $x^G$ in the power set $P(\Omega)$. 

One of the main tools in our study is 
the following ``Basic Lemma'' proved in Siemons and Volta \cite{SV}. The \emph{arity of a relation} $R$ is the set of cardinalities $ar(R) =\{|x| : x\in R\}$.

\bigskip\noindent \textbf{Lemma 3.1} \cite{SV} \hspace{1em} \textit{
Let  $H\leq S_n$ be a permutation group and $y$ a regular set in $H$. Suppose that one of the following holds:
\begin{enumerate}
\item[\rm (i)] There is a relation $R$ with $H = \G(R)$ so that $|y|  \notin ar(R)$, or
\item[\rm (ii)] 
$H$ is not set-transitive, and is maximal in $S_n$ with  this property. 
\end{enumerate}
Then
every subgroup of $H$ is a relation group.
}
\bigskip

Similarly as in \cite{SV}, we make essential use of the result by Seress \cite{ser} that describes all primitive permutation groups that have no regular set (cf.~\cite[Theorem~2.2]{SV}. There are only a finite number of them of degrees $5\leq n \leq 32$.

\bigskip\noindent \textbf{Theorem.} \hspace{1em} \textit{
Let $G$ be a simple primitive permutation group of degree $n$, other than the alternating group $A_n$. Then every subgroup $H$ of $G$ is a relation group with the following four exceptions:
\begin{enumerate}
    \item[\rm (i)] $n= 5$ and $G=H=C_5$,
    \item[\rm (ii)]  $n= 6$, $G=PSL(2,5)$, and $H=(C_5)^+$,
    \item[\rm (iii)] $n= 9$, $G=PSL(2,8)$, and  $H=PSL(2,8)$ or $H=AGL(1,8)^+$.
\end{enumerate}
}
\bigskip

The proof is based on extensive computation using system GAP 4.10.1. Below we give the theoretical part and some indications making possible performing the computations by the reader.

Since $G$ is by assumption primitive, by \cite[Theorem~4.2, Corollary~4.3]{SV} we know that if $G$ is other than $C_5$ or $PSL(2,8)$, then it is a relation group. 

Now, if $G$ is of degree $n \neq 5,6,9$, then by \cite{BP} (cf. \cite[Lemma~2.1]{SV}), $G$ is contained in a primitive group $G'$ that is maximal with the property of being not set-transitive. In turn, 
by \cite{ser} (cf. \cite[Theorem~2.2]{SV}), if $n$ is different from the cardinalities  $ n \leq 32$ listed in \cite[Theorem~2.2]{SV}, then $G'$ has a regular set. 
By Lemma~3.1 (ii) \cite{SV}, for the permutation group $G'$ having these two properties, every subgroup of $G'$ is a relational group.  Thus, it remains two check a finite number of simple permutation groups that appear in the lists in \cite[Theorem~2.2, Lemma~2.1]{SV}. 

In the next section we present the results of computations.
In this presentation we use a number of abbreviations. We write RG, meaning that a group in question is a \emph{relation group}. We say that a group $G$ is \emph{defined} by a set $x$ meaning that $G$ is defined by the orbit $x^G$ of $x$ in $G$, that is, $G =\G(x^G)$. We say that a group $G$ is defined in a group $G'$ by a set $x$ if $G =\G(R)$ with  $R=x^G\cup R'$ and $G'=\G(R')$. We also use some notation for the structure of permutation groups defined in \cite{kis}. Groups in question are assumed to act on a set $\Omega = \{1,2,\ldots,n\}$. Sets are denoted using square brackets. 

For each degree $n$ we list the simple primitive groups $G$ that (according to \cite[Theorem~2.2, Lemma~2.1]{SV}) have no regular set.
Then the computation is done recursively as follows. 
If $G$ is a relation group, then we point out a defining set, whose finding is based on computations of orbits of $G$ in $P(\Omega)$. 
If in addition $G$ has a regular set, satisfying one of the conditions of Lemma~3.1, then the computation terminates (since this means that each subgroup of the group in question is a relation group). 
Otherwise, we compute all (conjugate classes of) maximal subgroups of $G$ and repeat the whole procedure recursively for each such subgroup. Since the number of groups and subgroups is finite, and rather small, it is possible to complete the computation in a relatively short time. For most of the cases the depth of computation in the potential search tree does not exceed 2. There is only one case for the Mathieu group $M(12)$ of degree $12$, when one needs to inspect some subgroups as far as on the depth~$4$. 

The computations also show that the exceptions in the theorem not only fail to be relation groups, but they are not orbit closed in the sense of \cite{SV}. 

Our proof is based on ideas of the proof \cite[Theorem~4.1]{SV}. Yet, the reader should notice that there is a mistake in the proof of Theorem~4.1 in \cite{SV}. The argument works, only for $|\Omega|>32$, as in our proof.

Finally, note that in the theorem we consider only those simple permutation groups that are primitive (for which we have advanced tools based on the classification of finite simple groups). However, there are also imprimitive, and even intransitive actions of simple groups, for which we need other tools and more basic study to be done.

\section{Simple primitive groups with no regular set of degree $n$:}

\noindent $n=5$

There is only one group $G=C_5$. Not RG. No nontrivial subgroups. \bigskip

\noindent $n=6$

Only one group $G=PSL(2,5) = \langle (1,2,5)(3,4,6), (1,4)(2,5) \rangle$. Defined by $[1,2,5]$. Three maximal subgroups:

\begin{enumerate}

\item $G_{6,1}=\langle  (1,2,3)(4,5,6) ,  (1,5)(2,6)  \rangle$. Defined in $G$ by $[1,2]$.  Regular set $[3,4,5,6]$. All subgroups RG.

\item $G_{6,2}=\langle (1,2)(3,6), (1,6)(2,4) \rangle = D_5^+$. Defined in $G$ by $[1,2]$. No regular set.  Has two subgroups: $C_5^+$ not RG, and  $(S_2||S_2)^+$  RG.

\item $G_{6,3}=\langle (1,2,3)(4,5,6), (2,3)(5,6) \rangle = S_3 || S_3$.  Defined in $G$ by $[1,2]$. Regular set of size $4$.

\end{enumerate}
\smallskip

\noindent $n=7$

Only one group $G=PSL(3,2) = \langle (1,4)(6,7), (1,3,2)(4,7,5) \rangle$. Maximal not set-transitive. Defined by a set of size $2$. Three maximal subgroups:

\begin{enumerate}
\item $G_{7,1}=\langle  (2,3,4,7)(5,6), (2,5,3)(4,6,7)  \rangle$ transitive of order $24$. Defined in $G$ by a set of size $3$.  No regular set.

\begin{enumerate}
\item $H_{7,1} =\langle  (2,3,5)(4,7,6), (2,6,7)(3,4,5) \rangle$ of order $12$. 
Defined in $G_{7,1}$ by a set of size $4$. Regular set of
size~$5$.

\item $H_{7,2} =\langle  (2,4)(5,6), (2,4)(3,6,7,5) \rangle = (S_2 \oplus_\phi D_4)^+$ of order $8$. 
Defined in $G_{7,1}$ by a set of size $4$. 
Regular set of
size $5$.

\item $H_{7,3} =\langle  (2,7,5)(3,6,4), (3,6)(5,7) \rangle = (G_{6,3})^+$. 
\end{enumerate}
\smallskip

\item $G_{7,2}=\langle (2,4)(3,5,7,6), (1,2,4)(5,6,7)\rangle = S_3 \oplus_\phi S_4$ of order $24$.  
Defined in $G$ by a set of size $3$. No regular set.
Three maximal subgroups:

\begin{enumerate}
\item $H_{7,4} =\langle  (1,2,4)(3,5,7), (1,2,4)(3,7,6) \rangle  = C_3\oplus_\phi A_4$ of order $12$. 
Defined in $G_{7,2}$ by a set of size $4$. 
Regular set of size $5$.

\item $H_{7,5} =\langle  (3,6)(5,7), (2,4)(3,5,7,6) \rangle = (S_2\oplus_\phi D_4)^+$ of order $8$. 
Defined in $G_{7,2}$ by a set of size $4$. 
Regular set of size $5$.

\item $H_{7,6} =\langle  (1,4,2)(5,7,6), (2,4)(5,6) \rangle =  (G_{6,3})^+$. 
\end{enumerate}

\item $G_{7,3}=\langle (1,4,3)(2,6,7), (1,5,6)(2,7,3)\rangle $ of order $21$. Primitive RG. 
Two subgroups: $(C_3 || C_3 )^+$ and $C_7$; both RG. 
\end{enumerate}
\smallskip

\noindent $n=8$

Only one group 
$G=PSL(2,7) = \langle(1,6,5)(2,3,7), (1,4)(2,7)(3,5)(6,8)\rangle$. Defined by a set $x$ of size $4$ with $|x^G|=14$.
Three maximal subgroups:

\begin{enumerate}
\item  $G_{8,1} =\langle   (1,2,3,8)(4,6,5,7), (2,8,6)(3,4,5)\rangle $ of order $24$. Defined by a set of size $3$. Regular set of size $5$.

\item  $G_{8,2} =\langle    (1,6,5,2)(3,8,4,7) , (1,8,7)(2,6,4) \rangle $ of order 24. Conjugate in $S_8$ with  $G_{8,1}$. 

\item  $G_{8,3} =\langle  (2,7,5)(3,4,8), (2,8,4)(3,7,6) = (G_{7,3})^+ \rangle $ of order $21$. 
\end{enumerate}

\smallskip

\noindent $n=9$

Only one group 
$G=PSL(2,8) = \langle   (1,4,9)(2,7,5)(3,6,8), (1,3)(2,9)(5,7)(6,8)  \rangle$. Not RG. Three maximal subgroups:

\begin{enumerate}
\item $G_{9,1} = \langle  (1,8,3,9,5,2,4), (1,5)(2,6)(3,9)(4,8) \rangle = AGL(1,8)^+$, of order $56$. 
Not RG. The same orbits on subsets as $ASL(3,2)^+$. 
Two maximal subgroups:
\begin{enumerate}
\item $H_{9,1} = C_7\oplus I_2$. RG with no nontrivial subgroups.
\item $H_{9,2} = \langle   (1,2,3)(4,6,9)(5,7,8), (1,6)(2,4)(3,9)(7,8)  \rangle$. Defined by the union of orbits of  $[1,5],[1,2,3,7]$, and $[1,5,6,8,9]$. 
Regular set of size $3$.
\end{enumerate}

\item $G_{9,2}= \langle  (1,6)(2,4)(3,9)(7,8), (1,7)(2,5)(3,8)(4,6)  \rangle$ of order $18$. 
Defined by the union of orbits of  $[1,2,3,4]$ and $[1,2,3,4,5]$. 
Regular set of size $3$.

\item $G_{9,3}= \langle   (1,2)(3,5)(4,7)(6,9), (1,2)(3,8)(4,9)(5,6) \rangle$ of order  $14$. 
Defined by the union of orbits of  $[1,2,3,7]$ and $[1,5,6,8,9]$. 
Regular set of size $3$. 
\end{enumerate}

\bigskip\noindent $n=10$ 

Two groups: 

$G=(A_5,10) =  \langle  (1,3,5,7,9)(2,4,6,8,10),(1,9)(3,4)(5,10)(6,7)  \rangle $. Defined by sets $[1,2],[1,2,3]$. Regular set $[1,2,5,6]$. \\

$PSL(2,9)= \langle(2,5,9,7)(3,8,4,10), (1,3)(2,10)(4,7)(6,9)\rangle$ RG. No regular set. 
Five maximal subgroups:

\begin{enumerate}
\item $G_{10,1} \langle (1,7,5)(2,8,4)(3,10,9), (1,3)(2,10)(4,7)(6,9)\rangle =(A_5,10)$. Of order $60$. All as above. 
\item $G_{10,2} \langle (1,4,6)(2,9,8)(5,7,10), (1,3)(2,10)(4,7)(6,9)\rangle =(A_5,10)$. Of order $60$. As above.

\item $G_{10,3} \langle (1,4,5,2)(3,8,10,7), (1,4)(2,3)(5,10)(6,8)\rangle $. Of order $36$. Defined by $[1,2],[1,2,3,4]$. Regular set $[1,2,3]$.

\item $G_{10,4} \langle (3,9,6,10)(4,8,5,7), (1,10,3)(2,9,6)(5,8,7)\rangle $. of order 24.  Defined by $[1],[1,2],[1,2,3,4]$. Regular set $[5,\ldots,10]$.

\item $G_{10,5} \langle (1,8,2,7)(3,10,9,6), (1,4,7)(2,5,8)(3,6,10)\rangle $. of order 24. Conjugate in $S_{10}$ with $G_{10,4}$. 
\end{enumerate}

\bigskip\noindent $n=11$

Two groups: $PSL(2,11) \subseteq M_{11}$;
both without regular set. 

$M_{11}  \langle (1,10)(2,8)(3,11)(5,7), (1,4,7,6)(2,11,10,9)   \rangle$ defined by $[1,2,3,4,5]$ Five maximal subgroups.

\begin{enumerate}
\item $G_{11,1} \langle (1,5,6,2,7,11,4,9)(3,8), (1,7)(2,4)(6,9)(8,10)\rangle $. of order 48.    Defined by  the sets $[1],[1,2,3,4,5]$. Regular set $[1,2,3,4]$.

\item $G_{11,2} \langle(1,11,10,8,9,3)(2,6)(4,7,5), (1,10)(2,4)(5,6)(9,11)\rangle $. of order 120. Defined by the sets $[1], [1,2,3,4,5]$. Regular set $[ 1,2,3,7,8]$.

\item $G_{11,3} \langle (1,2,10,3,11,8)(4,5,9)(6,7), (1,11,5,4)(3,9,8,10)\rangle $. of order 144. The group of automorphisms of the set  $M(11)[1,2,3,4,5], G_{11,3}[1]$.  Regular set $[ 1,2,3,6]$.

\item $G_{11,4} \langle (3,10,6,11)(5,7,8,9), (2,10,3,7)(4,5,11,6)\rangle $. of order 720. 
The group of automorphisms of the set $M(11)[1,2,3,4,5], G_{11,3}[1]$.  No regular set. 
For maximal subgroups: 

\begin{itemize}

\item  $H_{11,1} =\langle   (2,4,8,5)(6,10,7,11), (2,10)(3,7)(5,6)(8,9) \rangle= PSL(2,9) \oplus I_1$. As for $n=10$.

\item $H_{11,2} =\langle  (3,4,11,8)(6,10,7,9), (2,6,9,3)(4,11,7,8) \rangle $ of order 72. 
The group of automorphisms of the set $M(11)[1,2,3,4,5], G_{11,3}[1],H_{11,2}[1,5]$. 
Regular set $[2,3,6]$.

\item $H_{11,3} = \langle  (3,11,6,10)(5,9,8,7), (2,10)(3,11)(4,7)(5,9)  \rangle $ of order 20. The group of automorphisms of the set $M(11)[1,2,3,4,5]$, $G_{11,3}[1]$, $H_{11,3} [12,10,3,11,6]$. Regular set $[2,4,5]$. 

\item $H_{11,4}  = \langle  (2,4,3,11)(5,10,6,7), (2,3)(4,5)(6,11)(8,9)  \rangle $. of order $16$. 
The group of automorphisms of the set  $M(11)[1,2,3,4,5]$, $G_{11,3}[1]$, $H_{11,4}[8,9])$. 
Regular set $[2,3,4,8]$.
\end{itemize}
\end{enumerate}

$PSL(2,11) \langle (2,3,4,10,5)(6,8,7,9,11), (1,5,4,8,11,9,3,6,10,2,7) \rangle$ defined by the set $[1,2,3,4,5]$.  
Four maximal subgroups: 

\begin{enumerate}
\item $G_{11,5}= \langle (1,5,2)(3,10,4)(8,11,9), (1,4)(2,9)(5,6)(8,10)  \rangle = (A_5,10)$. As above.

\item $G_{11,6}= \langle  (1,5,2)(3,10,4)(8,11,9), (1,5)(2,10)(6,11)(7,8) \rangle = A_5 \oplus_\phi  PSL(2,5)$. 
The group of automorphisms of the set $ PSL(2,11) [1,2,3,4,5]$, $G_{11,6} [1,2,6]$. 
Regular set $[1,2,6]$. 

\item $G_{11,7}= \langle  (1,5,8,10,2,11,3,9,6,7,4), (1,5,4,10,3)(6,9,7,11,8)\rangle$ of order $55$. 
The group of automorphisms of the set  $PSL(2,11)[1,2,3,4,5]$, $G_{11,7}  [1,24],$.  Regular set $[1,2]$.

\item $G_{11,8}= \langle
(1,9,6)(2,11)(3,5,4,7,8,10), (1,9)(3,10)(4,7)(5,8)
\rangle$ of order 12.  
The group of automorphisms of the set  $PSL(2,11)[1,2,3,4,5]$, $G_{11,7} [1,2,3]$. 
Regular set $[1,2,3]$. 
\end{enumerate}

\bigskip\noindent $n=12$

Three groups:
$M(12)$, $(M(11),12)$, $PSL(2,11)$.  First two have no regular set.

$M(12) = \langle (1,4,12,6)(2,7,5,9,8,10,3,11),  (1,12)(2,6,4,9,7,8,11,3)  \rangle$ Defined by the set $[1,2,3,4,5,6]$. 
11 maximal subgroups. $(M(11),11)^+$, $(M(11),12)$, $PSL(2,11)$ and:

\begin{enumerate}

\item $G_{12,1}= \langle
(1,5,11,8,6,7,4,3,9,2)(10,12), (2,8)(3,9)(6,7)(10,12)
\rangle$ of order $1440$.  
The group of automorphisms of the set  $M(12)[1,2,3,4,5,6]$, $G_{12,2} [1]$.  
No regular set. 6 maximal subgroups. 

\begin{itemize}

\item $H_{12,1}= \langle
(1,11,9,7)(2,3,5,6)(4,8)(10,12), (1,11,4,5)(2,9,8,3)
\rangle$ of order $720$.
Defined by three element set in $G_{12,2}$.
No regular set. 6 maximal subgroups: 

\begin{description}
 \item $F_{12,1} \langle
(1,4,11,5)(3,6,9,7), (1,5,8,6)(2,3,7,11)
\rangle = I_2 \oplus PSL(2,9)$, as above.

 \item $F_{12,2} \langle
(1,8,4,9)(2,3)(5,11,6,7)(10,12), (1,6)(3,11)(5,8)(10,12) 
\rangle$ of order $120$ 
Defined by two element set in $H_{12,1}$. 
Regular sets on 5 and 6 elements. 

 \item $F_{12,3} \langle
(1,8,4,9)(2,3)(5,11,6,7)(10,12), (1,2)(6,9)(8,11)(10,12)
\rangle$ Conjugate with $F_{11,2}$.

\item $F_{12,4} \langle
(1,5,3)(2,8,4,6,7,11)(10,12), (2,3)(5,6)(7,11)(10,12)
\rangle$ of order $72$.
Defined by two element set in $H_{12,1}$. 
Regular set $[1,2,3,10]$.

\item $F_{12,5} \langle
	(1,11,9)(2,4,6,5,8,3)(10,12), (1,9,7,11)(2,5)(3,4,6,8)(10,12)
\rangle$ of order $48$.
Defined by two element set in $H_{12,1}$. 
Regular set $[1,2,10]$.

\item $F_{12,6} \langle
(1,11,4,9,7,8)(2,5,3)(10,12), (1,11,9,7)(2,6,5,3)
\rangle$. 

Conjugate with  $F_{11,5}$.

\end{description}

\item $H_{12,2}= \langle
(1,2,9,3,4,7,6,8,11,5)(10,12), (1,6,9,3,11,8,2,7,4,5)(10,12)
\rangle$ of order $720$ 
Defined by four element set in  $G_{12,2}$. 
No regular set. 
4 maximal subgroups:

\begin{description}
 \item $F_{12,7} \langle
(2,6,7,8)(3,11,9,4), (1,3)(2,4)(5,7)(8,9)
\rangle = I_2 \oplus PSL(2,9)$, as above.

 \item $F_{12,8} \langle
(2,6,3,9,11,5,7,4)(10,12), (1,6)(2,4)(3,11)(7,9)
\rangle$ of order $72$.  
Defined by two element set in $H_{12,2}$. 
Regular set of size~3.
. 

 \item $F_{12,9} \langle
(1,3)(2,7)(4,8)(5,11)(6,9)(10,12), (1,9)(3,5)(4,11)(7,8)
\rangle$ of order $20$. 
Defined by two element set in $H_{12,2}$. 
Regular set of size~3.

\item $F_{12,10} \langle
(1,6,3,8,2,11,4,9)(10,12), (1,2)(5,7)(6,8)(9,11)
\rangle$ of order $16$.
Defined by two element set in $H_{12,2}$. 
Regular set of size 3.

\end{description}

\item $H_{12,3}= \langle
(2,4,3)(5,8,6)(7,9,11), (1,4,6,2)(3,7,11,9)
\rangle$ of order $720$. 
No regular set. 
Defined by five element set in $G_{12,1}$. 
For maximal subgroups:

\begin{description}
 \item $F_{12,11} \langle
(2,6,7,8)(3,11,9,4), (1,3)(2,4)(5,7)(8,9)
\rangle = I_2 \oplus PSL(2,9)$, as above.

 \item $F_{12,12} \langle
(2,3,11,7)(4,6,9,5), (1,6,9,3)(2,7,5,11)
\rangle$ of order $72$. 
Regular sets of size [3..9]. 
Defined by two element set in  $H_{12,3}$.  

 \item $F_{12,13} \langle
(2,4,11,9)(3,5,7,6), (1,9)(3,5)(4,11)(7,8)
\rangle$ of order $20$. 
Regular sets of size [2..10]. 
Defined by two element set in  $H_{12,3}$.  

\item $F_{12,14} \langle
(1,6,2,11)(3,9,4,8), (1,2)(5,7)(6,8)(9,11)
\rangle$ of order $16$. 
Regular sets of size [2..10]. 
Defined by two element set in $H_{12,3}$.  

\end{description}

\item $H_{12,4}= \langle
(1,6,8,2,9,11)(3,5,7)(10,12), (1,9,7,5)(2,3,11,8)
\rangle$ of size $144$.  
Defined by a three element set in $G_{12,2}$. 
Regular sets of size $[4..8]$.

\item $H_{12,5}= \langle
(1,5,8,4,11,3,6,2,7,9)(10,12), (1,5,7,2)(3,6,4,8)(9,11)(10,12)
\rangle$ of order $36$.   
Defined by three element set in $G_{12,2}$. 
Regular sets of size $[3..9]$.

\item $H_{12,6}= \langle
(1,9,3,6,2,8,4,11)(5,7), (3,4)(6,8)(9,11)(10,12), (3,4)(5,7)(6,9)(8,11)
\rangle$ of order $32$. 
Defined by three element set in $G_{12,2}$. 
Regular sets of size $[3..9]$.

\end{itemize}

\item $G_{12,2}= \langle
(1,5,11,8,6,7,4,3,9,2)(10,12), (1,9)(2,10)(3,7)(5,8)
\rangle$ of order $1440$.  The group of automorphisms of the sets $M(12)[1,2,3,4,5,6]$, $G_{12,3} [1,2,3,4,5]$.  
No regular set. 6 maximal subgroups.

\begin{itemize}
\item $H_{12,7}= \langle
(4,12,6,9)(5,8,10,7), (1,4)(2,8,7,10)(3,5)(6,9,12,11)
\rangle = S_6 \oplus_{\phi} S_6$  of order $720$. 
The group of automorphisms of the sets  $M(12)[1,2,3,4,5,6]$, $G_{12,2} [1,2,3,4,5]$, $H_{12,7}[1]$. 
No regular set. 

6 maximal subgroups: 

\begin{description}

\item $F_{12,15} \langle
(1,11)(2,7,3,8)(4,6,12,9)(5,10), (1,11,9,6)(2,10,8,5)(3,7)(4,12)
\rangle = A_6 \oplus_{\phi} A_6$. of order $360$. 
The group of automorphisms of the sets  $M(12)[1,2,3,4,5,6]$, $G_{12,2} [1,2,3,4,5],H_{12,7}[1]$, $F_{11,15}(X\setminus[1,2,3,4,9])$. No regular set.

5 maximal subgroups 

\begin{enumerate}

\item $P_{12,1} \langle (1,9,6)(4,11,12)(7,8,10), (1,6)(2,10)(3,8)(9,12)  \rangle$ of order $60$.
Regular sets of size $[3..10]$.

\item $P_{11,2} \langle (2,3,8)(5,10,7)(6,11,12), (1,6)(2,10)(3,8)(9,12) \rangle$. 
Conjugate to $P_{12,1} $.

\item $P_{12,3} \langle (1,6)(2,3,8,10)(4,12,11,9)(5,7), (2,5)(4,11)(6,9)(7,10)  \rangle$. of order $34$. 
Defined in $F_{12,15}$ by a set of size in $[3..10]$. 
Regular sets of size $[3..10]$.

\item $P_{12,4} \langle (1,4,11,6)(2,3)(5,10,7,8)(9,12), (1,12,4)(5,8,7)(6,11,9)  \rangle$. of order $24$. 
Defined in $F_{12,15}$ .
Regular sets of size $[3..10]$.

\item $P_{12,5} \langle (1,11,4,6)(2,5,3,7)(8,10)(9,12), (1,4,11)(2,8,5)(3,10,7) \rangle$. 
Conjugate to $P_{12,4} $. 
\end{enumerate}

\item $F_{12,16} \langle
(1,6,11,9)(2,3,8,7), (2,8)(3,7)(4,12)(5,10), (1,6)(2,8)(3,10)(5,7)
\rangle$ of order $48$.
Defined in $H_{12,7}$ by a set of any size. 
Regular sets of size $[3..10]$. 

\item $F_{12,17} \langle
(1,6,11,9,4)(2,5,3,10,8), (1,6)(2,8)(3,10)(5,7)
\rangle$ of order $120$.
Regular sets of size $[5..7]$. 
Defined in $H_{12,7}$ by a set of any size. 

\item $F_{12,18} \langle
(3,7,10,5)(4,9,12,11), (1,11,12,6,9,4)(2,8)(3,5,10)
\rangle$.
Conjugate to $F_{11,16} $.

\item $F_{12,19} \langle
(1,9)(4,6)(5,8)(11,12), (1,6)(2,3,7,8,10,5)(4,9,12)
\rangle$ of order $72$.
Regular sets of size $[4..8]$. 
Defined in $H_{12,7}$ by a set of size in $[2..11]$. 

\item $F_{12,20} \langle
(2,3,8,10)(4,9,11,12), (1,6,4)(3,7,8)(9,12,11)
\rangle$.
Conjugate to $F_{12,17} $.

\end{description}

\item $H_{12,8}= \langle
(1,5,11,8,6,7,4,3,9,2)(10,12), (1,4,9,12)(2,10,5,8)(3,7)(6,11)
\rangle$  of order $720$. 
Defined in $G_{12,2}$ by a set of size in $[4..9]$. 
No regular set.

4 maximal subgroups:

\begin{description}
\item  $F_{12,21} = \langle
(1,4,9,12)(2,10,5,8)(3,7)(6,11), (1,6)(2,10)(3,8)(9,12)
\rangle$ of order $360$. Conjugate to $F_{12,15}$. 

\item  $F_{12,22} = \langle
(1,2,12,10,11,5,6,3)(4,8,9,7), (2,5)(4,11)(6,9)(7,10)
\rangle$ of order $72$. 
Regular sets of size $[3..9]$. 
Defined in  $H_{12,8}$ by a set of size in $[2..10]$.

\item  $F_{12,23} = \langle
(1,2,6,7,12,5,11,3,4,10)(8,9), (2,10)(3,7)(4,6)(11,12)
\rangle$ of order $20$. 
Regular set of any size. 
Defined in $H_{12,8}$ by a set of any size.

\item  $F_{12,24} = \langle
(1,8,4,3,6,10,11,2)(5,9,7,12), (1,4)(2,3)(5,7)(6,11)
\rangle$ of order $16$. 
Regular set of any size 
Defined in $H_{12,8}$ by a set of any size.

\end{description}

\item $H_{12,9}= \langle
(1,4)(2,8,7,10)(3,5)(6,9,12,11), (1,5,6,2)(3,9)(4,10)(7,11,8,12)
\rangle$ of size $720$.
Defined in $G_{12,2}$ by a set of size in $[4..9]$. 
No regular set. 

4 maximal subgroups. 

\begin{description}

\item  $F_{12,25} = \langle
(1,4,9,12)(2,10,5,8)(3,7)(6,11), (1,6)(2,10)(3,8)(9,12)
\rangle$ Conjugate to $F_{12,15}$. 

\item  $F_{12,26} = \langle
(1,8)(2,11,5,4)(3,12)(6,7,9,10), (1,10,11,3)(2,6,5,12)(4,7)(8,9)
\rangle$. of size $72$. 
Regular sets of size $[3..9]$. 
Defined in $H_{12,9}$ by a set of size in $[2..10]$. 

\item  $F_{12,27} = \langle
(1,2,12,7)(3,6)(4,5,11,10)(8,9), (2,10)(3,7)(4,6)(11,12)
\rangle$ of order $20$. 
Regular set of any size. 
Defined in $H_{12,9}$ by a set of any size. 

\item  $F_{12,28} = \langle
(1,2,4,8,6,3,11,10)(5,9,7,12), (1,4)(2,3)(5,7)(6,11)
\rangle$ of order $16$.
Regular set of any size. 
Defined in $H_{12,9}$ by a set of any size.

\end{description}

\item $H_{12,10}= \langle
(1,6,4)(2,3,5,7,10,8)(11,12), (1,5,4,2)(3,9,7,12)(6,10)(8,11)
\rangle$ of order $144$.
Regular sets of size $[4..8]$.
Defined in $G_{12,2}$ by a set of size in $2..10$. 

\item $H_{12,11}= \langle
(1,2,11,8,4,7,12,5,6,10)(3,9), (2,5,10,8)(4,6,12,11)
\rangle$ of order $40$.
Regular sets of size $[3..9]$.  
Defined in $G_{12,3}$ by a set of any size.

\item $H_{12,12}= \langle
(1,3,12,7,9,2,11,10)(4,5,6,8), (1,9)(2,7)(3,10)(4,6), (1,9)$ $(2,10)(3,7)(5,8)
\rangle$  of order $32$.
Regular sets of size $[3..9]$. 
Defined in $G_{12,3}$ by a set of any size. 

\end{itemize}

\item $G_{12,3}= \langle
(1,2,10)(3,4,11)(5,12,9)(6,8,7), (1,10)(3,4)(5,8)(7,11)
\rangle$ of order $432$.  
The group of automorphisms of a set $G_{12,3} ([1,2,3,4,5,6],[1])$.
No regular set.
4 maximal subgroups. 

\begin{itemize}

\item $H_{12,13}= \langle
(4,7,12,9)(5,6,11,8), (1,2,10)(3,5,7)(4,8,6)(9,11,12)
\rangle$  of order $216$.
Regular sets of size $[4..8]$. 
Defined in $G_{12,3}$ by a set of size in $[5..7]$. 

\item $H_{12,14}= \langle
(4,6,12,8)(5,9,11,7), (2,10)(3,9)(4,5)(8,12)
\rangle$  of order $144$.
Regular sets of size $[4..8]$. 
Defined in $G_{12,3}$ by a set of any size. 

\item $H_{12,15}= \langle
(1,10,2)(4,12)(5,7,8,11,9,6), (2,10)(3,7)(4,6)(11,12)
\rangle$  of order $108$.
Regular sets of size $[4..8]$. 
Defined in $G_{12,3}$ by a set of size in $[2..10]$. 

\item $H_{12,16}= \langle
(1,10,2)(5,9,8)(6,11,7), (2,10)(4,11,7,8,12,5,9,6)
\rangle$  of order $48$.
Regular sets of sizes $[3..9]$. 
Defined in $G_{12,3}$ by a set of any size.

\end{itemize}

\item $G_{12,4}= \langle
(1,7,2)(3,6,11)(4,9,5)(8,12,10), (1,4)(2,12)(3,10)(7,11)
\rangle$ of order $432$.
The group of automorphisms of the set  $G_{12,4} ([1,2,3,4,5,6],[1,2,3,4])$.
No regular set. 
4 maximal subgroups. 

\begin{itemize}
\item $H_{12,17}= \langle
(1,11,2,4)(3,6)(5,9,8,10)(7,12), (1,4,5)(2,12,6)(3,10,9)(7,11,8)
\rangle$  of order $216$. 
Regular sets of size $[4..8]$. 
Defined in $G_{12,4}$ by a set of size $6$. \\
However, the automorphisms groups of the set $G_{12,4}[1,3],[1,2,3,4]$, $H_{12,17}[1,2,3,4,5,6]$. 

\item $H_{12,18}= \langle
(1,2,7)(4,6,11,5,12,8)(9,10), (1,11)(2,12,7,4)(3,8)(5,9,6,10)
\rangle$. of order $144$. 
Regular sets of size $[4..8]$. 
Defined in $G_{12,4}$  by a set of size in $[2..10]$. 

\item $H_{12,19}= \langle
(1,2,7)(4,6,11,5,12,8)(9,10), (2,7)(3,11,8)(4,5,9,12,6,10)
\rangle$. of order $108$. 
Regular sets of size $[4..8]$. 
Defined in $G_{12,4}$ by a set of any size. 

\item $H_{12,20}= \langle
(3,11,8)(4,6,9)(5,10,12), (1,11,3,8)(2,4,10,5,7,12,9,6)
\rangle$. of order $38$. 
Regular sets of size $[3..9]$. 
Defined in $G_{12,4}$ by a set of any size. 
\end{itemize}

\item $G_{12,5}= \langle
(1,10,12,6,4,3)(2,11,7,8,5,9), (1,11,4,9,8,12)(2,10,5,3,6,7)
\rangle$ of order $240$.
The group of automorphisms of the set  $G_{12,5} ([1,2,3,4,5,6],[1,2,3,4])$.
Regular set $[1,2,3,4,6]$.

\item $G_{12,6}= \langle
(1,5,12,3,6,4)(2,9,10)(8,11), (1,3)(2,7,9,10)(4,11,12,8)(5,6)
\rangle$ of order $192$. 
The group of automorphisms of the set $G_{12,6} ([1,3,4,5,6,8],[1])$.
Regular set $[6,7,8,10,11,12]$. It is the same size as one of the defining sets. However, since $[1]$ is one of the defining sets, they do not coincide.

\item $G_{12,7}= \langle
(1,4,7,3,6,10)(2,5,9,12,11,8), (1,6)(2,3)(4,9)(5,11)(7,10)(8,12)
\rangle$ of order $192$.
The group of automorphisms of the set $G_{12,7} ([1,2,3,4,5,6],[1,2,3,4,5])$.
Regular set $[4,7,8,9,10,11,12]$.

\item $G_{12,8}= \langle
(1,7,4,2,5,3)(6,9,8)(10,12), (1,9,2)(3,5,11)(4,8,10)(6,7,12)
\rangle$ of order $72$.
The group of automorphisms of the set  $G_{12,8} ([1,2,3,4,5,6],[1,2,3,4,5])$.
Regular set $[1,2,3,5]$. 

\end{enumerate}

$(M(11),12) \langle (1,12)(2,3,4,10)(5,9)(6,8,7,11), (1,9)(2,10)(3,7)(5,8)\rangle$.
The group of automorphisms of the set  $(M(11),12) ([1,2,3,4,6,10])$. 
5 maximal subgroups: $PSL(2,11)$ and

\begin{enumerate}

\item $G_{12,9}= \langle
(1,5)(2,8,6,7)(3,9)(4,10,11,12), (1,7)(2,10,8,4)(3,11)(5,6,12,9)
\rangle  = A_6 \wr S_2$ of order $720$. 
The group of automorphisms of the set  $G_{12,9} ([1,2,3,4,6,10]$ $,[1,2,3,4,5])$. 
No regular set. 
4 Maximal subgroups:

\begin{itemize}
\item $H_{12,21}= \langle
(1,12,10,11)(2,3,8,6)(4,5)(7,9), (1,11)(2,3)(5,10)(7,8)
\rangle$. of order $360$. Conjugate to $F_{12,15}$.

\item $H_{12,22}= \langle
(1,6)(2,10,7,11)(3,12)(4,9,5,8), (1,12,4,11)(2,6,7,9)(3,8)(5,10)
\rangle$. of order $72$. Conjugate to $F_{12,26}$.

\item $H_{12,23}= \langle
(1,3)(2,10)(4,9,12,7)(5,6,11,8), (1,11)(3,9)(4,5)(6,7)
\rangle$. of order $20$. Conjugate to $F_{12,27}$.

\item $H_{12,24}= \langle
(1,7,11,9,10,8,5,6)(2,4,3,12), (4,12)(5,11)(6,8)(7,9)
\rangle$. of order $16$. Conjugate to $F_{12,28}$.
\end{itemize}

\item $G_{12,10}= \langle
(1,11,9)(2,10,5,3,7,8)(4,6), (1,5)(2,11,7,9)(3,4)(6,8,12,10)
\rangle$ of order $144$. 
The group of automorphisms of the set  $G_{12,10} ([1,2,3,4,6,10],[1,2,3,4,5])$. 
Regular set $[1,2,3,5]$.

\item $G_{12,11}= \langle
(1,11,6)(2,7,4,8,3,12)(5,10), (1,9)(2,11)(3,12)(6,8)
\rangle$ of order $120$. y
The group of automorphisms of the set $G_{12,11} ([1,2,3,4,6,10],[1])$. 
Regular set $[1,2,3,4,5]$.

\item $G_{12,12}= \langle
(1,10,4,5,9,8,12,2)(3,6,7,11), (1,9)(4,8)(7,11)(10,12)
\rangle$ of order $48$. 
The group of automorphisms of the set  $G_{12,12} ([1,2,3,4,6,10],[1])$. 
Regular set $[1,2,3,7]$.

\end{enumerate}

$PSL(2,11) \langle (1,2,4)(3,9,12)(5,10,8)(6,7,11), (1,11)(2,8)(3,12)(4,9)(5,6)(7,10) \rangle$
The group of automorphisms of the set  $PSL(2,11)([1,2,3,4,5,6],[1,2,3,4])$. 
Regular set $[1,2,3,4,5]$

\bigskip\noindent $n=13$

Only one group: $PSL(3,3)= \langle (1,10,4)(6,9,7)(8,12,13), (1,3,2)(4,9,5)(7,8,12)(10,13,11)\rangle$  
The group of automorphisms of the set  $PSL(3,3)([1,2,3,4,5,6]$. 
4 maximal subgroups: 

\begin{enumerate}
\item $G_{13,1}= \langle
(1,11,9)(3,7,4)(5,6,10)(8,13,12), (2,12)(3,9)(4,7)(10,11)
\rangle$ of order $432$.  
The group of automorphisms of the set  $G_{13,1}([1,2,3,4,5,6],[1]$. 
No regular set. 
4 maximal subgroups.

\begin{itemize}
 \item $H_{13,1}= \langle (2,13)(3,10,9,4)(5,11,6,7)(8,12), (1,5,10)(2,8,12)(3,7,11)(4,6,9) \rangle$ of order 
 $216$. 
 Regular sets of size $[4..8]$. 
 Defined in $H_{13,1}$ by a set of size in $6,7$. 
 
  \item $H_{13,2}= \langle (2,8)(3,11,9,7)(4,6,10,5)(12,13), (1,4)(3,5)(7,9)(8,12)\rangle$ of order  
 $144$. 
 Regular sets of size $[4..8]$. 
 Defined in  $H_{13,1}$ by a set of size in $[2..10]$. 
 
  \item $H_{13,3}= \langle (2,12,8)(3,9)(4,11,5,10,7,6), (1,10)(3,11)(6,9)(8,12) \rangle$ of order 
 $108$. 
  Regular sets of size $[4..8]$. 
 Defined in $H_{13,1}$ by a set of any size. 
 
  \item $H_{13,4}= \langle (2,12,8)(4,7,5)(6,10,11), (2,8,13,12)(3,7,4,6,9,11,10,5) \rangle$. 
 $r = 48$. 
 Regular sets of size $[3..9]$. 
Defined in $H_{13,1}$ by a set of any size.

\end{itemize}

\item $G_{13,2}= \langle
(1,11,8)(3,13,10)(4,5,12)(6,7,9), (1,10)(6,12)(7,13)(8,9)
\rangle$ of order $432$. 
The group of automorphisms of the set   $G_{13,1}([1,2,3,4,5,6],[1]$. 
No regular set. 
4 maximal subgroups.

\begin{itemize}
 \item $H_{13,1}= \langle (1,4,10)(3,8,6)(5,12,7)(9,11,13), (1,11)(3,4,5,10)(6,13,9,8)(7,12) \rangle$ of order 
 $216$. 
 Regular sets of size $[4..8]$. 
Defined in $H_{13,1}$ by a set of size in $6,7$. 
 
  \item $H_{13,2}= \langle (3,5)(6,13)(7,12)(8,9), (1,9,4,6)(3,8)(5,13,11,12)(7,10)\rangle$ of order 
 $144$. 
 Regular sets of size$[4..8]$. 
 Defined in $H_{13,1}$ by a set of size in $[2..10]$. 
 
  \item $H_{13,3}= \langle (3,11)(6,8)(7,13)(9,12), (1,9,13)(4,7,12,10,6,8)(5,11) \rangle$ of order  
 $108$. 
  Regular sets of size $[4..8]$. 
 Defined in $H_{13,1}$ by a set of any size.
 
  \item $H_{13,4}= \langle (3,13,9)(5,8,6)(7,11,12), (1,13,3,9)(4,8,5,7,10,12,11,6) \rangle$ of order  
 $48$. 
 Regular sets of size $[3..9]$. 
Defined in  $H_{13,1}$ by a set of any size.

\end{itemize}

\item $G_{13,3}= \langle
(1,10,5)(2,12,8)(3,11,7)(4,9,6), (1,11,8)(3,13,10)(4,5,12)(6,7,9)
\rangle$ of order $39$. 
The groups of automorphisms of the set $G_{13,3}([1,2,3,4,5,6],[1,2]$. 
Regular set $X \setminus [1,2]$. 

\item $G_{13,4}= \langle
(1,11,4,5)(3,10)(6,12,9,13)(7,8), (1,5,11)(2,7,8)(3,12,6)(9,10,13)
\rangle$ of order $24$. 
The groups of automorphisms of the set $G_{13,3}([1,2,3,4,5,6],[1,2]$. 
Regular set $[1,2,6]$. 
\end{enumerate}

\bigskip\noindent $n=14$

Only one group: $PSL(2,13)= $ \\ $= \langle (1,7,14,5,11,9,3,6,13,2,12,4,8), (1,7,3,11,2,12,5,14,8,10,4,13,6) \rangle < PGL(2,13)$
Regular sets of size $6$ and $8$. 
The groups of automorphisms of the set $PGL(2,13)[1,2,3,4],$ 
$ PSL(2,13)[1,2,3]$. 

\bigskip\noindent $n=15$

Three groups:
$PSL(4,2)$, $(A(7),15)$ and $(A(6),15)$.

$PSL(4,2)=\langle (1,15,6,2,10)(3,5,7,13,12)(4,8,11,14,9),   (1,9)(2,5)(3,12)(4,11)(6,14)(10,13) $
Defining sets of size in $[3..12]$. 
No regular set.

6 maximal subgroups. $(A(7),15)$ and: 

\begin{enumerate}
 \item $G_{15,1}= \langle
(1,14,15)(2,6,12,9,13,7)(3,8)(4,10,5), (2,7)(3,6)(8,9)(10,14)(11,15)(12,13)
\rangle$ of order $1344$. 
No regular set. 
Defined in $PSL(4,2)$ by a set of any size. 
5 maximal subgroups: 

\begin{itemize}
 \item $H_{15,1} = \langle   (1,11,10)(2,8,9)(4,14,15)(7,13,12), (4,11)(5,10)(6,9)(7,8)  \rangle$ of order $168$. 
 Regular sets of size $[4..11]$. 
Defined in $G_{15,1}$ by a set of any size. 
 
  \item $H_{15,2} = \langle   (1,11,10)(2,9,3)(4,14,15)(6,7,12), (2,12)(3,13)(4,11)(5,10)(6,7)(8,9)  \rangle$ of order $168$. 
  Regular sets of size $[4..11]$. 
Defined in $G_{15,1}$ be a set of size in $[3..12]$.

   \item $H_{15,3} = \langle   (1,14,15)(2,6,3,8,12,9)(4,5,11)(7,13), (1,15)(2,7,13,6)(3,8,12,9)$
   $(4,5,10,11)  \rangle$ of order $192$. 
   Regular sets of size $[5..10]$. 
   Defined in $G_{15,1}$ by a set of any size.

    \item $H_{15,4} = \langle   	(1,4,14,11)(2,13)(3,9,12,6)(5,10), (1,10,14,5)(2,7,6,12)(3,13,8,9)$
    $(4,11)  \rangle$ of order $192$. 
    Regular sets of size $[5..10]$.
   Defined in $G_{15,1}$ by a set of any size.

     \item $H_{15,5} = \langle   (2,3)(4,11,14)(5,10,15)(6,8,12,7,9,13), (1,5,4)(3,7,6)(8,12,13)$
     $(10,14,15)  \rangle$. of order $168$. 
  Regular sets of size $[3..12]$
  Defined in  $G_{15,1}$ be a set of size in $[3..12]$. 
 
\end{itemize}

\item $G_{15,2}= \langle
(1,15)(2,5,7)(3,10,6,13,4,8)(9,12,11), (1,9)(2,4)(3,13)(5,11)(7,15)(10,12)
\rangle$ of order $1344$. 
No regular set. 
Defined in $PSL(4,2)$  by a set of any size.
5 maximal subgroups: 

\begin{itemize}
 \item $H_{15,6} = \langle   (2,6,5)(3,7,4)(8,11,12)(9,10,13), (1,3)(5,7)(9,11)(13,15) \rangle$ of order $168$. 
 Regular sets of size. $[4..11]$. 
 Defined in $G_{15,2}$ by a set of any size. 
 
  \item $H_{15,7} = \langle   (2,6,11)(3,7,10)(4,13,9)(5,12,8), (1,3)(4,10)(5,9)(6,8)(7,11)(13,15)  \rangle$ of order $168$. 
  Regular sets of size $[3..12]$. 
 Defined in $G_{15,2}$ by a set of size in $[3..12]$.

   \item $H_{15,8} = \langle   (2,7,13,8)(3,6,12,9)(4,11)(5,10), (1,15)(2,5)(3,10)(4,13)(6,8)(11,12)  \rangle$ of order $192$. 
   Regular sets of size $[5..10]$. 
   Defined in  $G_{15,2}$ by a set of any size.

    \item $H_{15,9} = \langle   	(1,7,15,9)(2,5,10,3)(4,13,12,11)(6,8), (1,8,9)(2,3,11)(5,12,13)$
    $(6,7,15)  \rangle$ of order $192$. 
    Regular sets of size $[5..10]$.
   Defined in  $G_{15,2}$ by a set of any size.

     \item $H_{15,10} = \langle  (1,5,6,15,11,8)(2,12)(3,9,4)(7,10,13), (1,13,6,9,10,2,5)$
     $(3,8,7,4,12,11,15) \rangle$ of order $168$. 
  Regular sets of size $[3..12]$
Defined in  $G_{15,2}$ by a set of any size $[3..12]$. 
 
\end{itemize}

\item $G_{15,3}= \langle
(1,15,6,2,10)(3,5,7,13,12)(4,8,11,14,9), (1,9)(2,5)(3,12)(4,11)(6,14)(10,13)
\rangle$ of order $720$. 
Regular sets of size $[6..9]$. 
Defined in $PSL(4,2)$ by a set of size in $[2..13]$.

\item $G_{15,4}= \langle
(1,6,4)(2,5,7)(8,9,15,11,10,12)(13,14), (1,13,6,10)(2,5)(3,8,4,15)(11,12)
\rangle$ of order $576$. 
Regular sets of size $[6..9]$. 
Defined in $PSL(4,2)$ by a set of any size.

\item $G_{15,5}= \langle
(1,2,11,6,7,5,14,8,15,10,4,12,3,9,13), (1,7)(2,12)(3,11)(4,10)(5,13)(9,15)
\rangle$ of order $360$. 
Regular sets of size $[4..11]$. 
Defined in $PSL(4,2)$ by a set of size in $[2..13]$.

\end{enumerate}

$(A(7),15)= \langle (1,15,6,2,10)(3,5,7,13,12)(4,8,11,14,9)$, \\ $(1,2,3)(4,9,13)(5,11,14)(6,10,12)(7,8,15) \rangle$. 

Regular sets of size $6$ and $8$. 
Defined in $PSL(4,2)$ by a set of size in $[4..11]$.

$(A(6),15) = \langle (1,14,12,3)(2,15)(4,8,11,10)(5,6,7,9)$, \\$(1,10)(2,9)(4,7)(5,13)(6,14)(12,15) \rangle$. It is a maximal subgroup of $(A(7),15)$. 

\bigskip\noindent $n=16$

There is no simple primitive group different than $A_{16}$.

\bigskip\noindent $n=17$

Only one group: $PSL(2,16) = \langle (3,13,9,10,11,7,17,5,4,12,14,6,8,16,15)$, \\ $(1,9,2)(3,14,6)(4,12,11)(7,16,10)(13,17,15) \rangle$. of order $4080$. 
Maximal subgroup of $P\Gamma L(2,16)$. 

Regular sets of size $5,7,8,9,10,12$.  \\
The groups of automorphisms of $P\Gamma L(2,16)[1,2,3,4,5]$, $PSL(2,16)[1,2,3,4]$, 

\noindent
$PSL(2,16)[1,2,3,4,5,6,13]$.

\bigskip\noindent $n=21$

Two groups:
$PSL(3,4) = M(21) =$  \\
$ \langle (1,7,12,16,19,21,6)(2,8,13,17,20,5,11)(3,9,14,18,4,10,15)$, \\ $(2,14,18,20,8)(3,7,12,13,19)(4,21,17,15,10)(5,11,16,6,9)\rangle$ of order $20160$ and

\noindent
$(A(7),21) =$ $\langle    (1,7,12,16,19,21,6)(2,8,13,17,20,5,11)(3,9,14,18,4,10,15)$, \\ $(4,6,5)(9,11,10)(13,15,14)(16,18,17)(19,20,21)   \rangle$.

$PSL(3,4)$ is a maximal subgroup \\ $P\Gamma L(3,4) = \langle   (1,5,6,19,17,14,4,10)(2,18,13,16)(3,12,21,15,20,8,7,9)$, \\  $(1,17,9,4,8,18,15,7,10,11,21,13,2,3,12)(6,20,16,14,19)   \rangle$ od order $120960$. 

The group of automorphisms of the set $P\Gamma L(3,4)[2,4,6,8,10],$  

\noindent $PSL(3,4)[1,2,3,4,5,6,13]$. 
Regular sets of the size $6,8,9,10,11,12,13,15$.

\noindent $(A(7),21)$ is a maximal subgroup\\ 
$(S(7),21)= \langle (1,7,12,16,19,21,6)(2,8,13,17,20,5,11)(3,9,14,18,4,10,15)$, \\ $(2,3)(4,5,6)(7,8)(9,10,11)(13,17,15,16,14,18)(19,21,20)   \rangle$ which is a maximal nor set transitive group and has a regular set.



\bigskip\noindent $n=22$

Only one group: $M(22) \langle (2,5,6,9,19,14,7)(3,21,17,18,4,13,22)$ , \\ $(1,20,15,2,13)(3,18,14,11,4)(5,6,10,8,7)(9,17,12,21,16)     \rangle$. 
This is a maximal subgroup of $M(22):2$, which is defined by a set of size in $[4..18]$. 
$M(22)$ is defined in  $M(22):2$ by a set of size in $[5..17]$. 

No regular set. 

8 maximal subgroups: 

\begin{enumerate}

\item    $G_{22,1} = \langle   (3,5,4,16)(6,12)(7,14,18,17)(8,19,21,15)(9,20,11,13)(10,22)$, \\ $(1,21,14)(3,17,9)(4,6,19)(5,8,15)(7,20,12)(11,22,16) \rangle$ of order 
$20160$.
This is $M(21) \oplus I_1$, see above.

\item    $G_{22,2} = \langle   (1,9,4,14,13)(2,12,21,7,18)(3,19,5,20,15)(6,17,11,10,16)$, \\ $(1,5)(2,20)(4,14)(6,17)(7,13)(8,16)(9,15)(19,22)  \rangle$ of order 
$5760$.

Regular sets of size $[5..17]$. 
Defined in $M(22)$ by a set of any size.

\item    $G_{22,3} = \langle (1,12)(2,18,15,21)(3,20,10,13)(5,7,22,17)(6,8,9,14)(11,19)$, \\ $(1,14)(2,18)(3,19)(4,9)(5,15)(6,17)(7,12)(11,16)   \rangle$ of order 
$2520$.

Regular sets of size $[5..17]$. 
Defined in $M(22)$ by a set of any size.

\item    $G_{22,4} = \langle (1,7,2,8,15,11)(3,6,5,4,14,10)(9,22,17)(12,21,16)(13,18)(19,20)$, \\ $(1,9)(3,22)(4,5)(6,17)(8,11)(12,13)(14,15)(18,20)   \rangle$ of order  
$2520$.
Conjugate to $G_{22,3}$. 

\item    $G_{22,5} = \langle  (1,8,11,2,21)(3,17,22,4,7)(5,20,19,6,18)(10,16,13,12,14)$, \\ $(1,5)(2,20)(4,14)(6,17)(7,13)(8,16)(9,15)(19,22)  \rangle$ of order 
$1920$.

Regular sets of size $[5..17]$. 
Defined in $M(22)$ by a set of any size.

\item    $G_{22,6} = \langle  (1,11,4,6)(2,13)(3,9)(5,14,19,22)(7,20,18,12)(8,15,16,21)$, \\ $(1,3)(2,9)(4,15)(5,7)(8,10)(11,17)(12,18)(13,21)   \rangle$ of order 
$1344$.

Regular sets of size $[4..18]$. 
Defined in $M(22)$ by a set of any size.

\item    $G_{22,7} = \langle  (1,21,18,5)(2,6,12,15,4,14,13,11)(3,10,8,20,7,9,17,16)(19,22)$, \\ $(1,21,16,8,9,17,10,7)(2,22,11,12,6,4,19,15)(3,18,5,20)(13,14)   \rangle$ of order  
$720$. 

Regular sets of size $[4..18]$. 
Defined in $M(22)$ by a set of any size.

\item    $G_{22,8} = \langle  (1,13,14,12,22,16,8,9,15,5,17)(2,18,4,10,11,21,7,19,20,3,6)$, \\ $(1,12)(2,21)(4,11)(5,9)(6,20)(7,10)(8,13)(14,22) \rangle$ of order 
$660$.

Regular sets of size $[4..18]$. 
Defined in $M(22)$ by a set of any size.
\end{enumerate}

\bigskip\noindent $n=23$

Only one group: \\
$M(23) = \langle (1,2,3,4,5,6,7,8,9,10,11,12,13,14,15,16,17,18,19,20,21,22,23)$, \\  $(2,16,9,6,8)(3,12,13,18,4)(7,17,10,11,22)(14,19,21,20,15)     \rangle$ of order  $10200960$. 
No regular set. 
Defining set of size in $[5..18]$. 

7 maximal subgroups:

\begin{enumerate}

\item    $G_{23,1} = \langle    (1,4,9,16)(2,3,17,20)(5,8,12,23)(6,14)(7,11,13,10)(15,19)$, \\ $(1,14,20)(2,19)(3,8,23,21,10,12)(4,15,17)(5,7,9,11,13,6)(16,18)  \rangle$ of order 
$443520 $. 
This is $M(22) \oplus I_1$, see above. 

\item    $G_{23,2} = \langle    (1,5,23,17,13,2,6,21)(3,7,10,11)(4,18)(9,14,16,12,20,22,19,15)$, \\ $(2,15)(4,18)(5,7)(8,13)(10,23)(11,12)(16,22)(17,19)  \rangle$ of order  
$40320  $. 
Regular sets of size $[7..16]$. 
Defined in $M(23)$ by a set of any size.

\item    $G_{23,3} = \langle    (1,4)(2,16,3,20,17,9)(6,7,13,8,19,15)(10,21,18)(11,22,12)(14,23)$, \\ $(2,15)(4,18)(5,7)(8,13)(10,23)(11,12)(16,22)(17,19)  \rangle$ of order 
$40320  $. 
Regular sets of size $[8..15]$. 
Defined in $M(23)$ by a set of any size.

\item    $G_{23,4} = \langle    (1,20)(3,7,8,22)(5,12)(6,18,17,19)(9,15,23,10)(11,16,13,14)$, \\ $(1,4)(2,9)(3,20)(6,13)(8,15)(10,21)(12,22)(16,17)  \rangle$ of order  
$20160 $. 
Regular sets of size $[7..16]$. 
Defined in $M(23)$ by a set of any size.

\item    $G_{23,5} = \langle    (1,5,21,18,4)(2,8,6,23,17)(3,7,10,16,12)(9,14,11,20,19)$, \\ $(1,5)(2,14)(4,16)(7,12)(8,19)(9,13)(10,15)(20,22)  \rangle$ of order 
$7920 $. 
Regular sets of size $[6..17]$. 
Defined in $M(23)$ by a set of any size.

\item    $G_{23,6} = \langle    (1,10,20,3,21,4)(2,9,19)(5,22,13,14,6,12)(7,18)(8,15)(16,17,23)$, \\ $(1,21)(3,23)(4,5)(6,16)(7,11)(8,14)(9,17)(10,19)  \rangle$ of order 
$5760 $. 
Regular sets of size $[6..17]$. 
Defined in $M(23)$ by a set of any size.

\item    $G_{23,7} = \langle    (1,15,22,8,14,18,12,2,19,3,17)(4,10,7,5,13,21,6,9,20,16,11)$, \\ $(1,21,12,20,11,7,23,14,22,15,16)(2,8,19,6,5,3,10,9,18,4,13)  \rangle$. 
$253$. 
Regular sets of size $[2..21]$. 
Defined in $M(23)$ by a set of size in  $[3..20]$.

\end{enumerate}

\bigskip\noindent $n=24$

Two groups: 
$M(24)$ i $PSL(2,23)$. The other one belongs to  $PGL(2,23)$, which is maximal not set transitive and has a regular set. 

$M(24)$ has a generating set of size in $[6..18]$.  
No regular set.

9 maximal subgroups:
\begin{enumerate}
\item $G_{24,1} = M(23) \oplus I_1$, see above.

\item $G_{24,2} = \langle    (1,21,3,7)(2,8,16,14)(4,12)(5,6)(9,10,22,24)(11,23,15,19)(13,20)$\\$(17,18),(1,18)(2,22)(4,12)(5,10)(7,16)(8,20)(13,19)(14,23) \rangle$ of order 887040. 

No regular set. 
Defined in $M(24)$ by a set of any size. 
7 maximal subgroups: 

\begin{itemize}

\item $H_{24,1} = \langle  (1,8,15,4)(3,6,13,9,7,11,17,19)(5,21,22,18,10,16,20,23)$, \\ $(14,24)(1,18,14,15,9,13,20)(3,8,5,6,11,17,22)(7,10,16,23,19,21,24)\rangle$ $= M(22) \oplus I_2$. of order $443520 $, see above.

\item $H_{24,2} = \langle  (1,22)(3,19,10,11)(4,24,17,7)(6,16,15,13)(8,20,18,23)(9,21)$, \\ $(2,12)(4,6)(7,19)(9,10)(13,14)(15,21)(17,22)(20,24)\rangle$ of order $40320 $. 

Regular sets of size $[7..17]$. 
Defined in $G_{24,2}$ by a set of any size.

\item $H_{24,3} = \langle  (1,7,17,6,14,18)(2,12)(3,24,4,5,10,16)(8,9,13)(11,21,20)$, \\ $(22,23)(1,10,4,7)(3,6,15,24)(5,16,8,14)(11,21)(13,18)(19,20,22,23)\rangle$ of order $11520 $. 

Regular sets of size  $[6..18]$. 
Defined in $G_{24,2}$ by a set of any size.

\item $H_{24,4} = \langle  (1,4,16,15,10,21,19,22,14,17)(2,12)(3,13,7,9,5,11,8,24,18,23)$, \\ $(6,20)(1,13)(3,5)(4,19)(6,20)(9,16)(10,14)(11,15)(17,23)\rangle$ of order $3840 $. 

Regular sets of size $[6..18]$. 
Defined in $G_{24,2}$ by a set of any size.

\item $H_{24,5} = \langle  (1,23,14,22)(2,12)(3,7,5,6)(4,9,24,13)(8,15)(10,11,16,21)$, \\ $(17,18)(19,20)(1,18)(4,19)(5,7)(6,11)(8,9)(14,24)(15,23)(17,20)\rangle$ of order $2688$. 
Regular sets of size $[5..19]$. 
Defined in $G_{24,2}$ by a set of any size.

\item $H_{24,6} = \langle  (1,23,10,7,19,3,24,20,9,11)(2,12)(4,5)$, \\ $(6,13,22,18,17,16,21,8,15,14)(2,12)(4,6)(7,19)(9,10)(13,14)(15,21)$ \\ $(17,22)(20,24)\rangle$ of order $1440$. 
Regular sets of size $[5..19]$. 
Defined in $G_{24,2}$ by a set of any size.

\item $H_{24,7} = \langle  (1,23,13,17)(2,12)(3,19,5,4)(6,10,20,14)(7,24)(8,21)$, \\ $(9,11,16,15)(18,22)(1,14)(3,5)(4,24)(6,7)(9,13)(10,16)(11,21)(22,23)\rangle$ of order $1320$. 

Regular sets of size $[4..20]$. 
Defined in $G_{24,2}$ by a set of size in $[2..22]$. 
\end{itemize}

\item $G_{24,3} = $ \\
$\langle  (1,19,15,14,24,21)(2,3,11,10)(4,6,22,5,20,8,18,12,7,13,17,23)(9,16)$, \\ $(1,14,9)(2,6,3)(4,13,22)(7,8,18)(10,12,17)(15,16,21) \rangle$ of order 
$322560 $

No regular set. 
Defined in  $M(24)$ by a set of any size. 
7 maximal subgroups:

\begin{itemize}

\item $H_{24,8} = \langle  (1,19,15,14,24,21)(2,7,17,11,22,20)(3,23,12)(4,18)(6,10,8)$, \\ $(9,16)(1,14,9)(2,6,3)(4,13,22)(7,8,18)(10,12,17)(15,16,21)\rangle $ of order $20160 $. 
Conjugate with $(G_{23,4})^+$

\item $H_{24,9} = \langle  (1,9)(2,23,8,11)(3,6,7,18)(4,22,12,10)(5,13,20,17)(14,19)$, \\ $(1,16,21)(2,5,17,8,6,12)(3,20,4,7,18,13)(9,24)(11,22,23)(15,19)\rangle$ of order $5760$. 
Regular sets of size $[6..18]$. 
Defined in $G_{24,2}$ by a set of any size.

\item $H_{24,10} = \langle  (1,9)(2,6,11)(3,22,23,20,5,8)(4,10,13,18,17,7)(14,16,21)$, \\ $(15,19)(1,9,24,16,14)(2,23,12,11,22)(3,18,4,13,7)(5,8,17,6,10)\rangle$ of order $11520 $. 
Regular sets of size $[6..18]$. 
Defined in $G_{24,2}$ by a set of any size.  

\item $H_{24,11} = \langle (1,16,21,24,15,14,9)(2,10,5,23,3,4,6,12,20,22,18,17,8,11)$, \\ $(7,13)(1,16,9,21,24)(2,10,4,18,7)(3,22,20,23,5)(6,11,12,17,13)\rangle$ of order $40320 $. 
Regular sets of size $[8..16]$. 
Defined in $G_{24,2}$ by a set of any size.

\item $H_{24,12} = \langle  (4,13)(5,18)(6,17)(7,22)(8,23)(9,21)(12,20)(14,19)$, \\ $(1,15,9,16,21,19)(2,4,3,11,12,10,23,6,22,8,18,7)(5,20,17,13)(14,24)\rangle$. of size $9216$. 
Regular sets of size $[6..18]$. 
Defined in $G_{24,2}$ by a set of size in $[2..22]$. 

\item $H_{24,13} = \langle (2,13,4,11,8,17,12,23)(3,10,6,5,7,22,18,20)(9,14,16,24)$, \\ $(19,21)(1,15,21,14,19,9,24)(2,22,13,23,20,6,4)(5,18,12,8,10,17,11)\rangle$ of size $21504$. 
Regular sets of size $[7..17]$. 
Defined in $G_{24,2}$ by a set of size in $[2..22]$. 

\item $H_{24,14} = \langle  (3,8,11,7,23,22,10)(4,18,20,13,6,17,12)(9,19,15,16,21,24,14)$, \\ $(1,16)(2,12,10,4,7,5,3,17,11,13,8,18)(6,23,20,22)(9,24,21,19,15,14)\rangle$ of order $21504$. 
Regular sets of size $[7..17]$. 
Defined in $G_{24,2}$ by a set of size in $[2..22]$. 

\end{itemize}

\item $G_{24,4} = $ \\
$\langle   (1,18,22,11,20,13,6,4)(3,17,10,19)(7,9,14,23,24,15,12,21)(8,16)$, \\ $(1,12)(2,17)(3,23)(4,21)(5,6)(7,13)(8,22)(9,19)(10,14)(11,16)(15,18)(20,24) \rangle$ of order 
$190080$

Regular sets of size $[10..14]$
Defined in $M(24)$ by a set of size in $[2..22]$.

\item $G_{24,5} = $ \\
$\langle   (1,24,5,15)(2,19,22,7,8,14,6,17,20,21,11,10)(3,4,12,13,23,18)(9,16)$, \\ $(1,16)(2,12)(3,6)(4,5)(7,23)(8,14)(9,21)(10,11)(13,15)(17,20)(18,24)(19,22)  \rangle$ of order 
$138240$.

Regular sets of size $[9..15]$
Defined in $M(24)$ by a set of size in $[2..22]$.

\item $G_{24,6} =$ \\
$ \langle   (1,21)(2,6,12,11,10,20,17)(4,19,22,15,9,8,23,7,24,18,14,16,13,5)$, \\ $(2,11,22)(3,21)(4,14,20,23,7,9)(5,19,10,15,24,13)(6,18)(8,16,12)  \rangle$ of order 
$120960$.

Regular sets of size $[9..15]$
Defined in $M(24)$ by a set of any size. 

\item $G_{24,7} = $ \\
$\langle   (1,13,24,14)(2,11)(3,17,15,7)(4,9,22,23)(5,19,12,6)(8,21)(10,16)(18,20)$, $(1,16)(2,12)(3,6)(4,5)(7,23)(8,14)(9,21)(10,11)(13,15)(17,20)(18,24)(19,22) \rangle$ of order 
$64512 $.

Regular sets of size $[7..17]$
Defined in $M(24)$ by a set of size in $[2..22]$.

\item $G_{24,8} = $ \\
$\langle  (1,13,12)(2,8,24)(3,10,15)(4,5,19)(6,18,9)(7,21,23)(11,14,22)(16,20,17)$, $(1,22,21,3)(2,19,12,5)(4,17,11,14)(6,20,23,16)(7,24,10,9)(8,15,18,13) \rangle$ of order  
$6072$. 

Regular sets of size $[5..19]$
Defined in $M(24)$ by a set of size in $[4..20]$.

\item $G_{24,9} = $ \\
$\langle   (1,11,18,6)(2,7,20,14)(3,19,17,23)(4,9,24,22)(5,8,15,21)(10,13,12,16)$, \\ $(1,15)(2,19)(3,9)(4,17)(5,11)(6,18)(7,16)(8,23)(10,13)(12,14)(20,21)(22,24)  \rangle$ of order  
$168$.

Regular sets of size $[3..21]$. 
Defined in $M(24)$ by a set of size in $[2..22]$. 

\end{enumerate}

\bigskip\noindent $n=32$

The only group is $PSL(2,31)$, which is a subgroup of maximal not-set-transitive group $PGL(2,31)$.
\bigskip





\nocite{*}
\bibliographystyle{aipnum-cp}%


\end{document}